\documentclass{amsart} 

\newtheorem*{mthm}{Main Theorem}
\newtheorem*{smthm}{Singh-Mehanna}
\newtheorem{thm}{Theorem}
\newtheorem{lemma}{Lemma}
\newtheorem*{cor}{Corollary}

\begin{document}

\title{Invariant Subspaces of $RL^1$}
\author{Daniel Jupiter}
\author{David Redett}
\address{Department of Mathematics, Texas A\&M University, College Station, TX 77843-3368}
\email{jupiter@math.tamu.edu}
\address{Department of Mathematics, Texas A\&M University, College Station, TX 77843-3368}
\email{redett@math.tamu.edu}
\thanks{The second author was supported in part by a VIGRE grant from the NSF}
\subjclass[2000]{Primary 47A15; Secondary 46E30}
\begin{abstract}
In this note we extend D. Singh and A. A. W. Mehanna's invariant subspace theorem for $RH^1$ (the real Banach space of analytic functions in $H^1$ with real Taylor coefficients) to the simply invariant subspaces of $RL^1$ (the real Banach space of functions in $L^1$ with real Fourier coefficients).
\end{abstract}
\maketitle

Let $\mathbf{T}$ denote the unit circle in the complex plane, and let $L^p$ denote the Lebesgue spaces on $\mathbf{T}$ with respect to Lebesgue measure normalized so that the Lebesgue measure of $\mathbf{T}$ is 1.  We use the standard notation $H^p$ to denote the subspace of $L^p$ consisting of those functions in $L^p$ whose negative Fourier coefficients vanish.  Let $$RH^p = \{ f \in H^p : \textnormal{the Fourier (Taylor) coefficients of f are real} \}.$$  An invariant subspace is a (closed) subspace invariant under multiplication by the coordinate function.  D. Singh and A. A. W. Mehanna \cite{S-M} gave a characterization of the invariant subspaces of $RH^1$.  Specifically, they proved the following result.

\begin{smthm}
Let $\mathcal{M}$ be an invariant subspace of $RH^1$.  Then there exists a unique (up to a constant multiple of modulus one) inner function, $I$, in $RH^1$  such that $\mathcal{M} = I\,RH^1$.
\end{smthm}

Let $$RL^p = \{ f \in L^p : \textnormal{the Fourier coefficients of f are real} \}.$$
A simply invariant subspace is an invariant subspace, $\mathcal{M}$, whose image under multiplication by the coordinate function is strictly contained in $\mathcal{M}$.  (In $RH^1$ every invariant subspace is simply invariant.)  In this note we extend Singh and Mehanna's result to the simply invariant subspaces of $RL^1$.

\begin{mthm}
Let $\mathcal{M}$ be a simply invariant subspace of $RL^1$.  Then there exists a unique (up to a constant multiple of modulus one) unimodular function, $U$, in $RL^1$  such that $\mathcal{M} = U\,RH^1$.
\end{mthm}

To prove this theorem we follow the approach of Singh and Mehanna.  We will first prove an analogous result in $RL^2$ and then use this result to prove the Main Theorem.  Singh and Mehanna's proof weighs heavily on the inner-outer factorization of functions in $H^1$.  In general, $L^1$ functions do not have such a factorization.  As we will soon see, however, the members of $\mathcal{M}$ have a nice factorization which will prove useful in the proof of the Main Theorem.

\begin{thm}\label{t:1}
Let $\mathcal{M}$ be a simply invariant subspace of $RL^2$.  Then there exists a unique (up to a constant multiple of modulus one) unimodular function, $U$, in $RL^2$  such that $\mathcal{M} = U\,RH^2$.
\end{thm}

To prove this, we need to understand the form of the doubly invariant subspaces of $RL^2$.  These are the invariant subspaces, $\mathcal{M}$, of $RL^2$ for which multiplication by the coordinate function takes $\mathcal{M}$ onto $\mathcal{M}$.

\begin{thm}
If $\mathcal{M}$ is a doubly invariant subspace of $RL^2$ then $\mathcal{M} = \mathbf{1}_{E}RL^2$ where $E$ is a measurable subset of $\mathbf{T}$.
\end{thm}

In the above theorem, $\mathbf{1}_{E}$ denotes the characteristic function of the set $E$.  The following proof is a slightly simplified version of the proof given in \cite{He} for the doubly invariant subspaces of $L^2$.  Although our proof is stated for $RL^2$, it works equally well in the $L^2$ setting.

\begin{proof}
If $1$ is in $\mathcal{M}$ then $\mathcal{M}=RL^2$.  In this case $E = \mathbf{T}$.  If $1$ is not in $\mathcal{M}$, let $q$ be the orthogonal projection of $1$ onto $\mathcal{M}$.  Then $1-q$ is orthogonal to $\mathcal{M}$ and, in particular, $1-q$ is orthogonal to $e^{in\theta}q$, for all integers $n$.  That is, $$\int^{\pi}_{-\pi}e^{in\theta}(q(e^{i\theta})-|q|^2(e^{i\theta}))\,d\theta = 0 \textnormal{\indent for all integers $n$.}$$  Therefore, $q=|q|^2$ (a.e.) on $\mathbf{T}$, and thus $q$ takes only values zero and one.  Let $E$ be the subset of $\mathbf{T}$ where $q$ takes the value one.  Then $q=\mathbf{1}_{E}$ and $1-q=\mathbf{1}_{E^c}$, so we have that $\mathbf{1}_{E} \in \mathcal{M}$ and $\mathbf{1}_{E^c} \in \mathcal{M}^{\perp}$.  By the double invariance of $\mathcal{M}$ we have that $\mathbf{1}_{E}RL^2$ is contained in $\mathcal{M}$ and that $\mathbf{1}_{E^c}RL^2$ is contained in $\mathcal{M}^{\perp}$.  We also have $\mathbf{1}_{E}RL^2 + \mathbf{1}_{E^c}RL^2 = RL^2$ and $\mathbf{1}_{E}RL^2 \cap \mathbf{1}_{E^c}RL^2 =\{0\}$.  Hence $\mathcal{M} = \mathbf{1}_{E}RL^2$ as desired.
\end{proof}

We are now ready to prove Theorem \ref{t:1}.  Our proof follows the proof given in \cite{He} for the simply invariant subspaces of $L^2$.  We include it for completeness.

\begin{proof}[Proof of Theorem 1]
Since $\mathcal{M}$ is simply invariant there exists a $U$ in $\mathcal{M} \ominus e^{i\theta}\mathcal{M}$ of norm one.  $U$ is orthogonal to $e^{in\theta}U$ for all natural numbers $n \geq1$.  So by the symmetry of the inner product on $RL^2$ we get $$\int_{-\pi}^{\pi}e^{in\theta}|U|^2(e^{i\theta})\,d\theta = 0 \textnormal{\indent for all integers $n \neq 0$}.$$  Thus $U$ has constant modulus one.  The set $\{e^{in\theta}U\}_{n=-\infty}^{\infty}$ spans a doubly invariant subspace in $RL^2$.  Since $U$ does not vanish on a set of positive measure, we have that this doubly invariant subspace is $RL^2$.  The span of $\{e^{in\theta}U\}_{n \geq 0}$ is $U\,RH^2$ and is contained in $\mathcal{M}$.  If we show that the set $\{e^{in\theta}U\}_{n < 0}$ is contained in $\mathcal{M}^{\perp}$, then we can conclude that $U\,RH^2$ is all of $\mathcal{M}$.  Showing that the set $\{e^{in\theta}U\}_{n<0}$ is contained in $\mathcal{M}^{\perp}$ is the same as showing that $U$ is orthogonal to $e^{in\theta}\mathcal{M}$ for all natural numbers $n>0$.  This is true by our choice of $U$.  Hence,  $\mathcal{M} = U\,RH^2$ as desired.  It remains to prove the uniqueness of $U$.  If $I$ is another unimodular function such that $\mathcal{M}=I\,RH^2$, then we have $U/I\,RH^2 = RH^2$.  Since the inverse of a unimodular function is its complex conjugate, we have that both $U\overline{I}$ and $\overline{U\overline{I}}$ are in $RH^2$.  This implies that $U/I$ is a constant of modulus one.
\end{proof}

Let $\mathcal{M}$ be a simply invariant subspace of $RL^1$.   Then $\mathcal{M}$ is a subset of $L^1$.  The complexification of  $\mathcal{M}$, $\overline{\mathcal{M}\otimes\mathbf{C}}^{L^1}$ is then a simply invariant subspace of $L^1$.  By a classical result \cite{He}, $\overline{\mathcal{M}\otimes\mathbf{C}}^{L^1}=\psi H^1$, where $\psi$ is a unimodular function in $L^1$.  So $\mathcal{M}$ is contained in $\psi H^1$.  It follows that every element of $\mathcal{M}$ has a unique unimodular-outer factorization.

Before we prove the Main Theorem we prove several technical lemmas.  Let $f$ be an element of $L^1$.  Define $f^*(e^{i\theta})=\overline{f(e^{-i\theta})}$.

\begin{lemma}\label{l:1}
For $f$ in $L^1$, $\overline{\hat{f}(n)}=\hat{f^*}(n)$.
\end{lemma}

\begin{proof}
Let  $f$ be an element of $L^1$.  Then 
\begin{eqnarray*}\hat{f^*}(n) &=& \int_{-\pi}^{\pi}f^*(e^{i\theta})e^{-in\theta}\,\frac{d\theta}{2\pi} \\
&=&  \int_{-\pi}^{\pi}\overline{f(e^{-i\theta})}e^{-in\theta}\,\frac{d\theta}{2\pi} \\
&=&  \overline{\int_{-\pi}^{\pi}f(e^{-i\theta})e^{in\theta}\,\frac{d\theta}{2\pi}} \\
&=& \overline{\int_{-\pi}^{\pi}f(e^{i\theta})e^{-in\theta}\,\frac{d\theta}{2\pi}} \\
&=& \overline{\hat{f}(n)}.
\end{eqnarray*}
\end{proof}

\begin{cor}
For $f$ in $L^1$, $f$ is in $RL^1$ if and only if $f=f^*$.
\end{cor}

\begin{lemma}
If $O$ is outer then $O^*$ is outer.
\end{lemma}

\begin{proof}
Recall that a function $f$ in $H^p$ is outer if and only if $$\log|f(0)|=\int_{-\pi}^{\pi}\log|f(e^{i\theta})|\,\frac{d\theta}{2\pi}.$$
Since $O$ is outer we have $$\log|O(0)|=\int_{-\pi}^{\pi}\log|O(e^{i\theta})|\,\frac{d\theta}{2\pi}.$$  By Lemma \ref{l:1}, $|O(0)|=|O^*(0)|$ and the negative Fourier coefficients of $O^*$ are zero.  A change of variable shows that $$\int_{-\pi}^{\pi}\log|O(e^{i\theta})|\,\frac{d\theta}{2\pi}=\int_{-\pi}^{\pi}\log|O^*(e^{i\theta})|\,\frac{d\theta}{2\pi}.$$  Thus $$\log|O^*(0)|=\int_{-\pi}^{\pi}\log|O^*(e^{i\theta})|\,\frac{d\theta}{2\pi},$$ so $O^*$ is outer.
\end{proof}

\begin{lemma}\label{l:2}
Let $f$ be an element of $RL^1$ such that $f$ has a unique factorization $f=UO$, where $U$ is a unimodular function and $O$ is an outer function.  Then $U$ and $O$ are in $RL^1$.
\end{lemma}

\begin{proof}
Since $f$ is in $RL^1$, by the above corollary we have that $f=f^*$.  Thus $UO=U^*O^*$.  Since $U^*$ is unimodular and $O^*$ is outer, the uniqueness of our factorization implies $U^*=U$ and $O^*=O$.  By the above corollary we get $U$ and $O$ are in $RL^1$.
\end{proof}

We are now ready to prove the Main Theorem.

\begin{proof}[Proof of Main Theorem]
Let $\mathcal{M}$ be a simply invariant subspace of $RL^1$.  Then $\mathcal{M} \cap RL^2$ is an invariant (closed) subspace of $RL^2$.  We will show that $\mathcal{M} \cap RL^2$ is dense in $\mathcal{M}$.  Then $\mathcal{M} \cap RL^2$ is actually simply invariant, and by Theorem \ref{t:1} $\mathcal{M} \cap RL^2=U\,RH^2$.  Hence $\mathcal{M}$ is of the form $U\,RH^1$, as desired.

We first show that $\mathcal{M} \cap RL^2$ is nonempty.  Let $f$ be any nonzero element of $\mathcal{M}$.  By Lemma \ref{l:2} we have that $f=UO$, where $U$ is unimodular and $O$ is outer, with both $U$ and $O$ in $RL^1$.  Since $O$ is outer it is actually a member of $RH^1$,  so by Lemma 3.4 of \cite{S-M} we may assume without loss of generality that $\sqrt{O}$ is in $RH^2$.   Therefore $g := U \sqrt{O}$ is in $RL^2$.  We now show that $g$ is also in $\mathcal{M}$.  By Corollary 3.4 of \cite{S-M} there exists a sequence of polynomials, $\{p_{n}\}$, in $RH^2$ such that $$\|\sqrt{O}p_{n}-1\|_{2} \rightarrow 0 \textnormal{ as $n\rightarrow \infty$}.$$  Thus, 
\begin{eqnarray*}
\|fp_{n} - g \|_{1} & = & \|g\sqrt{O}p_{n}-g\|_{1} \\
& \leq & \|g\|_{2}\,\|\sqrt{O}p_{n}-1\|_{2} \textnormal{\indent (by Cauchy-Schwarz)} \\
&\rightarrow & 0 \textnormal{ as $n \rightarrow \infty$}.
\end{eqnarray*}  Since $fp_{n}$ is in $\mathcal{M}$ for all $n$ by the invariance of $\mathcal{M}$, and since $\mathcal{M}$ is closed, we see that $g$ is in $\mathcal{M}$, as desired.

It remains to show that $\mathcal{M} \cap RL^2$ is dense in $\mathcal{M}$.  Let $f$ be any nonzero element of $\mathcal{M}$.  Then $f=UO$, where $U$ is unimodular and $O$ is outer, with both $U$ and $O$ in $RL^1$.  As mentioned above, we assume without loss of generality that $\sqrt{O}$ is in $RH^2$.  Let $\sqrt{O}_{n} = \sum_{k=0}^{n}a_{k}e^{ik\theta}$ be the partial sums of the Fourier series for $\sqrt{O}$.  We know that $\sqrt{O}_{n}$ converges to $\sqrt{O}$ in $RL^2$.  By the work above we know that $U\sqrt{O}$ is in $\mathcal{M} \cap RL^2$.  By the invariance of $\mathcal{M} \cap RL^2$ we have that $U\sqrt{O}\sqrt{O}_{n}$ is in $\mathcal{M} \cap RL^2$ for all $n>0$.  Thus
\begin{eqnarray*}
\|U\sqrt{O}\sqrt{O}_{n} - f\|_{1} &=& \|U\sqrt{O}(\sqrt{O}_{n}-\sqrt{O})\|_{1} \\
& \leq & \|U\sqrt{O}\|_{2}\|\sqrt{O}_{n}-\sqrt{O}\|_{2} \textnormal{ \indent (by Cauchy-Schwarz)} \\
& \rightarrow & 0 \textnormal{ as $n \rightarrow \infty$}.
\end{eqnarray*}
We conclude that $\mathcal{M} \cap RL^2$ is dense in $\mathcal{M}$, as desired.
\end{proof}

\end{document}